\documentclass{article}%
\usepackage{amsmath}
\usepackage{amsfonts}
\usepackage{amssymb}
\usepackage{graphicx}%
\begin{document}

\begin{center}
\textbf{\ Boundedness on inhomogeneous Lipschitz spaces of fractional
integrals, singular integrals and hypersingular integrals associated to
non-doubling measures.}

\bigskip

\ A. \ Eduardo Gatto\footnote{This paper was supported by a sabbatical leave
from DePaul University and by a grant from the Ministerio de Educacion y
Ciencia of Spain, Ref. \# SAB2006-0118 for a quarter sabbatical visit to the
Centre de Recerca Matem\`{a}tica in Barcelona.
\par
The author wants to express a deep gratitude to these institutions, to the
Analysis group of the Universidad Aut\'{o}noma de Barcelona, in particular to
Xavier Tolsa for the invitation and to the director and staff of the CRM for
their pleasant hospitality and excellent resources.\pagebreak
}

DePaul University\medskip\ and CRM

\ \ \ \ \ \ \ \ \ \ \ \ \ \ \ \ \ \ \ \ \ \ \ \ \ \ \ \ \ \ \ \ \ \ \ \ \ \ \ \ \ \ \ \ \ \ \ \ \ \ \ \ \ \ \ \ \ \ \ \ \ \ \ \ \ \ \ \ \ \ \ $To$
$my$ $sons$
\end{center}

\bigskip\ \ \ \ \ \ \ \ \ \ \ \ \ \ \ \ \ \ \ \ \ \ \ \ \ \ \ \ \ \ \ \ \ \ \ \ \ \ \ \ \ \ \ \ \ \ \ \ \ \ \ \ \ \ \ \ \ \ \ \ \ \ \ \ \ \ \ \ \ \ \ \ \ \ \ \ \ \ \ \ \ \ \ \ \ \ \ \ \ \ \ \ \ \ \ \ \ \ \ \ \ \ \ \ \ \ \ \ \ \ \ \ \ \ \ \ \ \ \ \ 

\ \ \ \ \ \ \ \ \ \ \ \ \ \ \ \ \ \ \ \ \ \ \ \ \ \ \ \ \ \ \ \ \ \ \ \ \ \ \ \ \ \ \ \ \ \ \ \ \ \ \ \ \ \ \ \ \ \ \ \ \ \ \ \ \ \ \ \ \ \ \ \ \ \ \ \ \ \ \ \ \ \ \ \ \ \ \ \ \ \ \ \ \ \ \ \ \ \ \ \ \ \ \ \ \ \ \ \ \ \ \ \ \ \ \ \ \ \ \ \ \ \ \ \ \ \ \ \ \ \ \ 

\bigskip

\textbf{Abstract}:

In the context of a finite measure metric space whose measure\ satisfies a
growth condition, we prove $"T1"$ type necessary and sufficient conditions for
the boundedness of fractional integrals, singular integrals, and hypersingular
integrals on inhomogeneous Lipschitz spaces. We also indicate how the results
can be extended to the case of infinite measure. Finally we show applications
to Real and Complex Analysis. \bigskip\ \bigskip\bigskip\bigskip
\vspace{6.05cm}\medskip

\newpage

\noindent1. \textbf{Introduction. Definitions and Statement of the Theorems}

\qquad Let (X,d,$\mu$) be a finite measure metric space whose measure $\mu$
satisfies a n- dimensional growth condition, that is, (X,d) is a metric space
and $\mu$ is a finite Borel measure that satisfies the following condition:
there is $n>0$ and a constant $A>0$ such that $\mu(B_{r})\leq Ar^{n},$for all
balls $B_{r}$ of radius $r$ and for all $r>0$. Note that this condition allows
the consideration of non-doubling as well as doubling measures.

\qquad Our results will apply to functions defined on the support of $\mu,$ of
course the support of $\mu$ has to be well defined, where supp$(\mu)$ is the
smallest closed set $F$ such that for all Borel sets $E$, $E\subset F^{c}%
,\mu(E)=0.$ For example, if $X$ is separable, then the support of $\mu$ is
well defined. Furthermore to avoid any confusion we will assume that $X= $
$\operatorname*{supp}(\mu)$

\qquad The inhomogeneous Lipschitz-H\"{o}lder spaces of order $\beta
,0<\beta\leq1,$ will be denoted $\Lambda_{\beta}$ and consists of all bounded
functions $f$ \ that satisfy

$\sup_{x\neq y\in X}\frac{\left\vert f(x)-f(y)\right\vert }{d^{\beta}%
(x,y)}<\infty.$ \ The space $\Lambda_{\beta}$ is a Banach space with the norm
$\left\Vert f\right\Vert _{\Lambda_{\beta}}=\sup_{x\in X}\left\vert
f(x)\right\vert +\sup_{x\neq y\in X}\frac{\left\vert f(x)-f(y)\right\vert
}{d^{\beta}(x,y)}.$ \ It will be useful to have a notation for each term in
the norm, let $\sup(f)=\sup_{x\in X}\left\vert f(x)\right\vert $, and
$\left\vert f\right\vert _{\beta}=\sup_{x\neq y\in X}\frac{\left\vert
f(x)-f(y)\right\vert }{d^{\beta}(x,y)}.$

\qquad The results in this paper have extensions to the case $\mu(X)=\infty$,
but the constants depend on the normalization of the integrals at infinity, we
will indicate these extensions after the section on proofs. The letter $C$,
$c$ will denote constants not necessarily the same at each ocurrence.\bigskip

Let $\Omega=X\times X\backslash\Delta,$ where $\Delta=\left\{
(x,y):x=y\right\}  .$ \ A function $L_{\alpha}(x,y):\Omega\rightarrow$C will
be called a standard fractional integral kernel of order $\alpha,0<\alpha
<1,$when there are constants $B_{1}$ and $B_{2}$ such that

\begin{enumerate}
\item[(L1)] $\left\vert L_{\alpha}(x,y)\right\vert \leq\frac{B_{1}%
}{d^{n-\alpha}(x,y)}$.

\item[(L2)] $\ \left\vert L_{\alpha}(x_{1},y)-L_{\alpha}(x_{2},y)\right\vert
\leq B_{2}\frac{d^{\gamma}(x_{1},x_{2})}{d^{n-\alpha+\gamma}(x_{1},y)}$ , for
some $\gamma,\alpha<\gamma\leq1,$ and $2d(x_{1},x_{2})\leq d(x_{1},y).$
\end{enumerate}

The fractional integral of order $\alpha$ of a function $f$ in $\Lambda
_{\beta}$ is defined by:%

\[
L_{\alpha}f(x)=\int L_{\alpha}(x,y)f(y)d\mu(y).
\]
Note that in particular $L_{\alpha}(x,y)=\frac{1}{d^{n-\alpha}(x,y)}$ is a
standard fractional kernel of order $\alpha$.

\bigskip

\noindent\textbf{Theorem 1}

Let $0<\alpha<\gamma\leq1,$ $0<\beta<1,$ and $\alpha+\beta\leq1$ when $1<n$ or
$\alpha+\beta<n$ when $n\leq1$. \ The following statements are equivalent:

\begin{enumerate}
\item[a)] $L_{\alpha}1\in\Lambda_{\alpha+\beta}$.

\item[b)] $L_{\alpha}:\Lambda_{\beta}\rightarrow\Lambda_{\alpha+\beta} $ is bounded.
\end{enumerate}

We define now the singular integral kernels that we will consider in Theorem 2
and Theorem 3. \ A function $K(x,y):\Omega\rightarrow C$ will be called a
standard singular integral kernel when there are constants $C_{1,}C_{2}$ and a
number $\gamma,0<\gamma\leq1,$ such that

\begin{enumerate}
\item[(S1)] $\left\vert K(x,y)\right\vert \leq\frac{C_{1}}{d^{n}(x,y)}$

\item[(S2)] $\left\vert K(x_{1,}y)-K(x_{2},y)\right\vert \leq C_{2}%
\frac{d^{\gamma}(x_{1},x_{2})}{d^{n+\gamma}(x_{1},y)},$ for $2d(x_{1}%
,x_{2})\leq d(x_{1},y)$
\end{enumerate}

Let $\eta$ be a function in $C^{1}\left[  0,\infty\right)  $ such that
$0\leq\eta\leq1,\eta(s)=0$ for $0\leq s\leq1/2$ and $\eta(s)=1$ for $1\leq s.$
\ Let $K_{\varepsilon}(x,y)=\eta(\frac{d(x,y)}{\varepsilon})K(x,y),$
$\varepsilon>0$ where $K(x,y)$ is a standard singular integral kernel . We
will denote $T_{\varepsilon}$ the operator $T_{\varepsilon}f(x)=\int
K_{\varepsilon}(x,y)f(y)d\mu(y).$

\noindent\newline\textbf{Theorem 2}

Let $K(x,y)$ be a standard singular integral kernel. \ Let $0<\beta
<\min(n,\gamma).$ \ The following two statements are equivalent:

\begin{enumerate}
\item[a)] $\left\Vert T_{\varepsilon}1\right\Vert _{\Lambda_{\beta}}\leq C$,
for all $\epsilon>0.$

\item[b)] $T_{\varepsilon}:\Lambda_{\beta}\rightarrow\Lambda_{\beta}$ are
bounded and $\left\Vert T_{\varepsilon}\right\Vert _{\Lambda_{\beta
}\rightarrow\Lambda_{\beta}}\leq C^{\prime},$ for all $\varepsilon>0.$
\end{enumerate}

One of the novelties in this Theorem is that the cancellation condition (S3)
for all $x$ (see below) follows from part b).

\bigskip

In Theorem 3 we will consider Principal Value Singular Integrals. \ We will
denote by $Lip_{\beta}$ the space of classes of measurable functions $f$ for
which there is a $g\in\Lambda_{\beta}$ such that $f=g$ $\ $except for a set
$E$ that depends on $f$, with $\mu(E)=0$. The norm of $f$ in $Lip_{\beta}$ is
defined as $\left\Vert f\right\Vert _{Lip_{\beta}}=\left\Vert f\right\Vert
_{\infty}+\left\vert f\right\vert _{\beta}$, where $\left\vert f\right\vert
_{\beta\text{ \ }}=\sup_{x\neq y\in X}\frac{\left\vert g(x)-g(y)\right\vert
}{d^{\beta}(x,y)}=\sup_{x\neq y\in X-E}\frac{\left\vert f(x)-f(y)\right\vert
}{d^{\beta}(x,y)}.$

We also need to add the following two conditions on the kernel:

\begin{enumerate}
\item[(S3)] \bigskip$\left\vert \int_{r_{1}<d(x,y)<r_{2}}K(x,y)d\mu
(y)\right\vert \leq C_{3}$ for all $0<r_{1}<r_{2}<\infty,\mu-a.e$ in $x.$

\item[(S4)] $\lim_{\varepsilon\rightarrow0}\int_{\varepsilon<d(x,y)<1}%
K(x,y)d\mu(y)$ \ exists $\mu-a.e$ in $x.$
\end{enumerate}

The principal value singular integral of a function $f\in Lip_{\beta}$ is
defined by%

\[
Kf(x)=\lim_{\varepsilon\rightarrow0}\int_{\varepsilon<d(x,y)}K(x,y)f(y)d\mu
(y)
\]

\noindent\textbf{Theorem 3}

\qquad Let $K(x,y)$ be a standard singular integral kernel that in addition
satisfies (S3) and (S4). \ Let $0<\beta<\min(n,\gamma)$ and $f\in Lip_{\beta
}.$ \ Then $Kf(x)$ is well defined $\mu-a.e.$ and the following two statements
are equivalent:

\begin{enumerate}
\item[a)] $K1\in Lip_{\beta}$

\item[b)] $K:Lip_{\beta}\rightarrow Lip_{\beta}$ is bounded
\end{enumerate}

A function $D_{\alpha}(x,y):\Omega\rightarrow C$ will be called a standard
hypersingular kernel of order $\alpha,0<\alpha<1,$when there are constants
$E_{1}$ and $E_{2}$ such that:

\begin{enumerate}
\item[(D1)] $\left\vert D_{\alpha}(x,y)\right\vert \leq\frac{E_{1}%
}{d^{n+\alpha}(x,y)},$

\item[(D2)] $\left\vert D_{\alpha}(x_{1},y)-D_{\alpha}(x_{2},y)\right\vert
\leq E_{2}\frac{d^{\gamma}(x_{1},x_{2})}{d^{n+\alpha+\gamma}(x_{1},y)},$ for
some $\gamma,0<\gamma\leq1,$ and $2d(x_{1},x_{2})\leq d(x_{1},y).$
\end{enumerate}

\bigskip

The hypersingular integral of order $\alpha$ of a function $f\in\Lambda
_{\beta}$ $\alpha<\beta\leq1$is defined by:
\[
D^{\alpha}f(x)=\int D_{\alpha}(x,y)\left[  f(y)-f(x)\right]  d\mu(y)
\]
Note that in particular $D_{\alpha}(x,y)=\frac{1}{d^{n+\alpha}(x,y)}$ is a
standard hypersingular kernel of order $\alpha$ when $X=R^{n}$ and $\mu$ is
the Lebesgue measure, and we have $\int\frac{1}{d^{n+\alpha}(x,y)}\left[
f(y)-f(x)\right]  dy=c_{\alpha}(\Delta^{\frac{\alpha}{2}}f)(x)$ for $f$
sufficiently smooth and $0<\alpha<2.\left[  \text{S}\right]  $

\bigskip

\noindent\textbf{Theorem 4}

\qquad Let $0<\alpha<\beta\leq1$ and $\beta-\alpha<n.$Then $D^{\alpha}$:
$\Lambda_{\beta}\rightarrow\Lambda_{\beta-\alpha}$ is bounded.

\bigskip

Note that $D^{\alpha}1=0$ by definition. Also, Theorem 4 and its proof are
valid without changes in the case $\mu(X)=\infty.$

\bigskip

\bigskip

\noindent\textbf{2. Proofs}\newline\newline

We would like to point out that the proofs are based on classical methods, see
for example [Z], adjusted to the modern "T1" formulation and to the present
general context. \ For carrying out the proofs we need the following known
lemma about measures that satisfy the n-dimensional growth condition.

\bigskip

\noindent\textbf{Lemma}

Let $(X,d,\mu)$ be a measure metric space such that $\mu$ satisfies the
n-dimensional growth condition, $r>0.$ Then

\begin{enumerate}
\item $\int_{d(x,y)<r}\frac{1}{d^{n-\delta}(x,y)}d\mu(y)\leq c_{1}r^{\delta},$
\ $0<\delta\,<n.$

\item $\int_{r\leq d(x,y)}\frac{1}{d^{n+\delta}(x,y)}d\mu(y)\leq
c_{2}r^{-\delta},$ \ \ $0<\delta$

\item $\int_{r/2\leq d(x,y)<r}\frac{1}{d^{n}(x,y)}d\mu(y)\leq c_{3}$
\end{enumerate}

\bigskip

\noindent\textbf{Proof of the Lemma}

The three parts are a consequence of the growth condition. To prove part1, we
rewrite the integral as a series and mayorize each term using the growth
condition and we add the resulting series. \ In detail we have:%

\begin{align*}
\int_{d(x,x_{o})<r}\frac{1}{d^{n-\delta}(x,x_{o})}d\mu(x)  &  =\sum
_{k=0}^{\infty}\int_{2^{-k-1}r\leq d(x,x_{o})<2^{-k}r}\frac{1}{d^{n-\delta
}(x,x_{o})}d\mu(x)\leq\\
\sum_{k=0}^{\infty}\frac{\mu(B_{2^{-k-1}r}(x_{o}))}{(2^{-k-1}r)^{n-\delta}}
&  \leq A\sum_{k=0}^{\infty}\frac{(2^{-k}r)^{n}}{(2^{-k-1}r)^{n-\delta}%
}=Ar^{\delta}(\frac{2^{n}}{2^{\delta}-1}).
\end{align*}

To prove part 2 we perform a similar estimate:%

\begin{align*}
\int_{d(x,x_{o})\geq r}\frac{1}{d^{n+\delta}(x,x_{o})}d\mu(x)  &  =\sum
_{k=0}^{\infty}\int_{2^{k}r\leq d(x,x_{o})<2^{k+1}r}\frac{1}{d^{n+\delta
}(x,x_{o})}d\mu(x)\leq\\
\sum_{k=0}^{\infty}\frac{\mu(B_{2^{k+1}r}(x_{o}))}{(2^{k}r)^{n+\delta}}  &
\leq A\sum_{k=0}^{\infty}\frac{(2^{k+1}r)^{n}}{(2^{k}r)^{n+\delta}%
}=Ar^{-\delta}(\frac{2^{n}2^{\delta}}{2^{\delta}-1})
\end{align*}

Finally for part 3 we have:%

\[
\int_{r/2\leq d(x,y)<r}\frac{1}{d^{n}(x,y)}d\mu(y)\leq\frac{\mu(B_{r}(x_{o}%
))}{(r/2)^{n}}\leq A2^{n}.
\]

\bigskip

\noindent\textbf{Proof of Theorem 1}

Observe first that $1\in\Lambda_{\beta}$ and therefore condition b) implies
condition a). \ We will prove now that condition a) implies condition b). \ We
can just consider the case $L_{\alpha}(x,y)$=$\frac{1}{d^{n-\alpha}(x,y)}%
,$because the general case is proven in the same way, and we will denote
$L_{\alpha}=I_{\alpha}$.

Condition (L1) is clearly valid. To show that condition (L2) is verified ,we
use the Mean Value Theorem Consider $2d(x_{1},x_{2})\leq d(x_{1},y)$, and
$0<\theta<1$ we have:%

\begin{align*}
\left\vert \frac{1}{d^{n-\alpha}(x_{1},y)}-\frac{1}{d^{n-\alpha}(x_{2}%
,y)}\right\vert  &  \leq\sup_{\theta}\left\vert (-n+\alpha)(\theta
d(x_{1},y)+(1-\theta)(d(x_{2},y))^{-n+\alpha-1}\right\vert \text{ .}\\
\left\vert d(x_{1},y)-d(x_{2},y)\right\vert  &  \leq B_{2}\frac{d(x_{1}%
,x_{2})}{d^{n-\alpha+1}(x_{1},y)}%
\end{align*}

Now we will estimate $\sup(I_{\alpha}f).$ \ Let $x\in X.$ We will use the
Lemma to obtain%

\[
\left\vert I_{\alpha}f(x)\right\vert \leq\int_{d(x,y)<1}\frac{\left\vert
f(y)\right\vert }{d^{n-\alpha}(x,y)}d\mu(y)+\int_{1\leq d(x,y)}\frac
{\left\vert f(y)\right\vert }{d^{n-\alpha}(x,y)}d\mu(y)\leq\left\Vert
f\right\Vert _{\infty}(c_{1}+\mu(X)),
\]
and therefore $\sup(I_{\alpha}f)\leq\sup(f)(c_{1}+\mu(X)).$ \ We will estimate
next $\left\vert I_{\alpha}f\right\vert _{(\beta)}.$ We write%

\[
I_{\alpha}f(x_{1})-I_{\alpha}f(x_{2})=\int_{X}\frac{f(y)}{d^{n-\alpha}%
(x_{1},y)}d\mu(y)-\int_{X}\frac{f(y)}{d^{n-\alpha}(x_{2},y)}d\mu(y)=\int
_{X}\frac{f(y)-f(x_{1})}{d^{n-\alpha}(x_{1},y)}d\mu(y)+
\]

\begin{align*}
&  f(x_{1})\int_{X}\frac{1}{d^{n-\alpha}(x_{1},y)}d\mu(y)-\int_{X}%
\frac{f(y)-f(x_{1})}{d^{n-\alpha}(x_{2},y)}d\mu(y)-f(x_{1})\int_{X}\frac
{1}{d^{n-\alpha}(x_{2},y)}d\mu(y)=\\
&  \int_{X}\frac{f(y)-f(x_{1})}{d^{n-\alpha}(x_{1},y)}d\mu(y)-\int_{X}%
\frac{f(y)-f(x_{1})}{d^{n-\alpha}(x_{2},y)}d\mu(y)+f(x_{1})\left[  I_{\alpha
}1(x_{1})-I_{\alpha}1(x_{2})\right]  .
\end{align*}
The last term can be mayorized using the hypothesis, and we have $\left\vert
f(x_{1})\left[  I_{\alpha}1(x_{1})-I_{\alpha}1(x_{2})\right]  \right\vert \leq
c\sup(f)d^{\alpha+\beta}(x_{1},x_{2}).$

Let now $r=d(x_{1},x_{2})$ and $B_{2r}(x_{1})$ the ball of radius $2r$ and
center $x_{1}.$ We write%

\[
\left\vert \int_{X}\frac{f(y)-f(x_{1})}{d^{n-\alpha}(x_{1},y)}d\mu(y)-\int
_{X}\frac{f(y)-f(x_{1})}{d^{n-\alpha}(x_{2},y)}d\mu(y)\right\vert \leq
\]

\[
\int_{B_{2r}(x_{1})}\frac{\left\vert f(y)-f(x_{1})\right\vert }{d^{n-\alpha
}(x_{1},y)}d\mu(y)+\int_{B_{2r}(x_{1})}\frac{\left\vert f(y)-f(x_{1}%
)\right\vert }{d^{n-\alpha}(x_{2},y)}d\mu(y)+
\]

\[
\int_{B_{2r}^{c}(x_{1})}\left\vert f(y)-f(x_{1})\right\vert \left\vert
\frac{1}{d^{n-\alpha}(x_{1},y)}-\frac{1}{d^{n-\alpha}(x_{2},y)}\right\vert
d\mu(y)=J_{1}+J_{2}+J_{3}
\]
For the first term using the lemma we have%

\[
J_{1}\leq\left\vert f\right\vert _{(\beta)}\int_{B_{2r}(x_{1})}\frac{d^{\beta
}(x_{1},y)}{d^{n-\alpha}(x_{1},y)}d\mu(y)\leq c\left\vert f\right\vert
_{(\beta)}r^{\alpha+\beta}=c\left\vert f\right\vert _{(\beta)}d^{\alpha+\beta
}(x_{1},x_{2}).
\]
For the second term we write%

\[
J_{2}\leq\left\vert f\right\vert _{(\beta)}\int_{B_{3r}(x_{2})}\frac
{2r}{d^{n-\alpha}(x_{2},y)}d\mu(y)\leq c\left\vert f\right\vert _{(\beta
)}d^{\alpha+\beta}(x_{1},x_{2}),
\]
For the third term we use (L2) and the lemma to get%

\[
J_{3}\leq\left\vert f\right\vert _{(\beta)}\int_{B_{2r}^{c}(x_{1})}\frac
{B_{2}}{d^{n-\alpha-\beta}(x_{1},y)}d\mu(y)\leq c\left\vert f\right\vert
_{(\beta)}d^{\alpha+\beta}(x_{1},x_{2})
\]
Collecting the previous estimates, we have%

\[
\left\Vert I_{\alpha}f\right\Vert _{\Lambda_{\beta}}\leq C\left\Vert
f\right\Vert _{\Lambda_{\beta}}
\]
This concludes the proof of Theorem 1.

\bigskip

\bigskip

\noindent\textbf{Proof of Theorem 2}

Observe first that $1\in\Lambda_{\beta}$ and therefore condition b) implies
condition a).

Before doing the proof of the theorem and for the sake of completeness, we
will show that $K_{\varepsilon}$ satisfies conditions (S1) and (S2) with
constants independent of $\varepsilon.$

Condition (S1) is true because $\eta$ is bounded. \ To show condition (S2),
assume that $2d(x_{1},x_{2})\leq d(x_{1},y)$ and consider the following two cases:

\begin{enumerate}
\item[Case 1:] $1<\frac{d(x_{1},y)}{\varepsilon}$ and $1<\frac{d(x_{2}%
,y)}{\varepsilon}.$ In this case $K_{\varepsilon}(x,y)=K(x,y),$ and therefore
(S2) is true with the same constant.

\item[Case 2:] $1\geq\frac{d(x_{1},y)}{\varepsilon}$ or $1\geq\frac
{d(x_{2},y)}{\varepsilon}.$ \ Assume $1>\frac{d(x_{1},y)}{\varepsilon}.$ \ 
\end{enumerate}

\noindent We write%

\begin{align*}
\left\vert K_{\varepsilon}(x_{1},y)-K_{\varepsilon}(x_{2},y)\right\vert  &
\leq\left\vert \eta(\frac{d(x_{1},y)}{\varepsilon})-\eta(\frac{d(x_{2}%
,y)}{\varepsilon})\right\vert \left\vert K(x_{1},y)\right\vert +\\
&  \left\vert \eta(\frac{d(x_{2},y)}{\varepsilon})\right\vert \left\vert
K(x_{1},y)-K(x_{2},y)\right\vert
\end{align*}
The first term above is less than or equal to%

\begin{align*}
\left\Vert \eta^{\prime}\right\Vert _{\infty}\frac{\left\vert d(x_{1}%
,y)-d(x_{2},y)\right\vert }{\varepsilon}\left\vert K(x_{1},y)\right\vert  &
\leq\left\Vert \eta^{\prime}\right\Vert _{\infty}\frac{d(x_{1},x_{2}%
)}{\varepsilon}\left\vert K(x_{1},y)\right\vert \leq\\
c(\frac{d(x_{1},x_{2})}{\varepsilon})^{\gamma}\left\vert K(x_{1}%
,y)\right\vert  &  \leq c\frac{d^{\gamma}(x_{1},x_{2})}{d^{n+\gamma}(x_{1},y)}%
\end{align*}
On the other hand the second term is less than or equal to $c\left\vert
K(x_{1},y)-K(x_{2},y)\right\vert \leq c\frac{d^{\gamma}(x_{1},x_{2}%
)}{d^{n+\gamma}(x_{1},y)}.$ \ If $1\geq\frac{d(x_{2},y)}{\varepsilon}$ the
proof is similar.

\bigskip

\bigskip To show that condition a) implies condition b), the first step\ is to
get the cancellation (S3) of the kernel, for all $x\in X$.

\bigskip Observe that for $0<r_{1}<r_{2}<\infty$, we have

\qquad\qquad%
\begin{align*}
\qquad T_{r_{1}}1(x)-T_{r_{2}}1(x)  &  =\int_{\frac{1}{2}r_{1}<d(x,y)\leq
r_{1}}\eta(\frac{d(x,y)}{r_{1}})K(x,y)d\mu(y)+\int_{r_{1}<d(x,y)}%
K(x,y)d\mu(y)\\
&  -\int_{\frac{1}{2}r_{2}<d(x,y)<r_{2}}\eta(\frac{d(x,y)}{r_{2}}%
)K(x,y)d\mu(y)-\int_{r_{2}\leq d(x,y)}K(x,y)d\mu(y)
\end{align*}
Since the left hand side is uniformly bounded in $r$ and $x$, and also the
first and third terms are uniformly bounded because of the growth condition
(see lemma), it follows that

\bigskip{}%

\[
(S3)\ \ \ \ \ \ \ \ \ \ \ \ \ \ \ \ \ \ \ \ \ \ \ \ \ \ \left\vert \int
_{r_{1}<d(x,y)<r_{2}}K(x,y)d\mu(y)\right\vert \leq C,\text{ for all }x.
\]
\ \ \ 

\bigskip

\bigskip Now, we will estimate $\sup\left\vert T_{\in}f(x)\right\vert
.$Observe first that%

\[
T_{\varepsilon}f(x)=\int_{d(x,y)\leq1}K_{\epsilon}(x,y)f(y)d\mu(y)+\int
_{d(x,y)>1}K_{\epsilon}(x,y)f(y)d\mu(y)=
\]

\bigskip%

\[
\int_{d(x,y)\leq1}K_{\epsilon}(x,y)f(y)-f(x)d\mu(y)+f(x)\int_{\frac{1}%
{2}\varepsilon<d(x,y)\leq\varepsilon}K_{\epsilon}(x,y)d\mu(y)
\]

\qquad\qquad\qquad%

\[
f(x)\int_{\varepsilon<d(x,y)\leq1}K(x,y)d\mu(y)+\int_{d(x,y)>1}K_{\epsilon
}(x,y)f(y)d\mu(y)
\]

\qquad\qquad\qquad

Now, by conditions (S1), (S3) and the growth condition we can bound the
absolute value\ of the terms\ above by $\left\Vert f\right\Vert _{\Lambda
_{\beta}}$ and therefore $\sup_{x\in X}\left\vert T_{\varepsilon
}f(x)\right\vert \leq c\left\Vert f\right\Vert _{\Lambda_{\beta}}.$

\qquad\qquad\qquad

\qquad\qquad\qquad

Next, \ we estimate $\sup_{x\neq y}\frac{\left\vert T_{\varepsilon
}f(x)-T_{\varepsilon}f(y)\right\vert }{d^{\beta}(x,y)}.$ \ We consider the
difference $T_{\varepsilon}f(x_{1})-T_{\varepsilon}f(x_{2})$, and the
following decomposition:\bigskip%

\[
T_{\varepsilon}f(x_{1})-T_{\varepsilon}f(x_{2})=\int K_{\varepsilon}%
(x_{1},y)f(y)d\mu(y)-\int K_{\varepsilon}(x_{2},y)f(y)d\mu(y)
\]
\[
=\int K_{\varepsilon}(x_{1},y)\left[  f(y)-f(x_{1})\right]  d\mu
(y)+f(x_{1})\int K_{\varepsilon}(x_{1},y)d\mu(y)-
\]

\[
\int K_{\varepsilon}(x_{2},y)\left[  f(y)-f(x_{1})\right]  d\mu(y)-f(x_{1}%
)\int K_{\varepsilon}(x_{2},y)d\mu(y)=
\]

\[
\int K_{\varepsilon}(x_{1},y)\left[  f(y)-f(x_{1})\right]  d\mu(y)+\int
K_{\varepsilon}(x_{2},y)\left[  f(y)-f(x_{1})\right]  d\mu(y)+
\]

\[
f(x_{1})\left[  T_{\varepsilon}1(x_{1})-T_{\varepsilon}1(x_{2})\right]
\]
\qquad\qquad

Observe now that the last term can be estimated using the hypothesis and we have

$\qquad\qquad\qquad\qquad\qquad$%
\[
\left\vert f(x_{1})\left[  T_{\varepsilon}1(x_{1})-T_{\varepsilon}%
1(x_{2})\right]  \right\vert \leq c\sup(f)d^{\beta}(x_{1},x_{2}).
\]

To estimate the first two terms, let $r=d(x_{1},x_{2}),$ we rewrite them as follows:%

\[
\int K_{\varepsilon}(x_{1},y)\left[  f(y)-f(x_{1})\right]  d\mu(y)+\int
K_{\varepsilon}(x_{2},y)\left[  f(y)-f(x_{1})\right]  d\mu(y)=
\]

\[
\int_{d(x_{1},y)<3r}K_{\varepsilon}(x_{1},y)\left[  f(y)-f(x_{1})\right]
d\mu(y)+\int_{d(x_{1},y)<3r}K_{\varepsilon}(x_{2},y)\left[  f(y)-f(x_{1}%
)\right]  d\mu(y)+
\]

\[
\int_{3r<d(x_{1},y)}\left[  f(y)-f(x_{1})\right]  \left[  K_{\varepsilon
}(x_{1},y)-K_{\varepsilon}(x_{2},y)\right]  d\mu(y)=H_{1}+H_{2}+H_{3}
\]

\qquad\qquad

The absolute value of $H_{3}$ can be estimated as follows,%

\[
\left\vert H_{3}\right\vert \leq\left\vert f\right\vert _{\beta}d^{\gamma
}(x_{1},x_{2})\int_{3r<d(x_{1},y)}\frac{d^{\beta}(x_{1},y)}{d^{n+\gamma}%
(x_{1},y)}d\mu(y)\leq c\left\vert f\right\vert _{\beta}d^{\beta}(x_{1},x_{2})
\]
For $\left\vert H_{1}\right\vert $we have%

\[
\left\vert H_{1}\right\vert \leq\left\vert f\right\vert _{\left[
\beta\right]  }\int_{d(x_{1},y)<3r}\frac{C_{1}}{d^{n-\beta}(x_{1},y)}%
d\mu(y)\leq c\left\vert f\right\vert _{\left[  \beta\right]  }d^{\beta}%
(x_{1},x_{2})
\]
Finally to estimate $H_{2}$ we write%

\[
\int_{d(x_{1},y)<3r}K_{\varepsilon}(x_{2},y)\left[  f(y)-f(x_{1})\right]
d\mu(y)=
\]

\[
\int_{d(x_{1},y)<3r}K_{\varepsilon}(x_{2},y)\left[  f(y)-f(x_{2})\right]
d\mu(y)
\]

\[
+\left[  f(x_{2})-f(x_{1})\right]  \int_{\left\{  y:\varepsilon/2<d(x_{2}%
,y)\right\}  \cap\left\{  y:d(x_{1},y)<3r\right\}  }K_{\varepsilon}%
(x_{2},y)d\mu(y)=J_{1}+J_{2}
\]
For the first term we have%

\[
\left\vert J_{1}\right\vert \leq\int_{d(x_{2},y)<4r}\frac{c\left\Vert
f\right\Vert _{(\beta)}}{d^{n-\beta}(x_{2},y)}d\mu(y)\leq c\left\vert
f\right\vert _{\beta}d^{\beta}(x_{1},x_{2})
\]
To estimate $J_{2}$ consider first%

\[
\int_{d(x_{1},y)<3r}K_{\varepsilon}(x_{2},y)d\mu(y)=\int_{d(x_{2}%
,y)<2r}K_{\varepsilon}(x_{2},y)d\mu(y)+
\]

\[
\int_{\left\{  y:d(x_{1},y)<3r\right\}  \backslash\left\{  y:d(x_{2}%
,y)<2r\right\}  }K_{\varepsilon}(x_{2},y)d\mu(y)
\]
Observe now that%

\[
\left\vert \int_{d(x_{2},y)<2r}K_{\varepsilon}(x_{2},y)d\mu(y)\right\vert \leq
C_{3}
\]
\bigskip and using part 3 of the lemma we get%

\[
\left\vert \int_{\left\{  y:d(x_{1},y)<3r\right\}  \backslash\left\{
y:d(x_{2},y)<2r\right\}  }K_{\varepsilon}(x_{2},y)d\mu(y)\right\vert \leq
\int_{\left\{  y:2r<d(x_{2},y)<4r\right\}  }\left\vert K_{\varepsilon}%
(x_{2},y)\right\vert d\mu(y)\leq c
\]
therefore%

\[
\left\vert J_{2}\right\vert \leq c\left\vert f\right\vert _{\beta}d^{\beta
}(x_{1},x_{2})
\]
collecting the estimates we have:%

\[
\left\vert K_{\varepsilon}f(x_{1})-K_{\varepsilon}f(x_{2})\right\vert \leq
c\left\Vert f\right\Vert _{\Lambda_{\beta}}d^{\beta}(x_{1},x_{2})
\]
and finally%

\[
\left\Vert K_{\varepsilon}f\right\Vert _{\Lambda_{\beta}}\leq c\left\Vert
f\right\Vert _{\Lambda_{\beta}},
\]
with c independent of $\varepsilon.$

\bigskip

\bigskip

\noindent\textbf{Proof of Theorem 3}

Observe first that $1\in$ $Lip_{\beta}$ and therefore condition b) implies
condition a). \ Let $f\in Lip_{\beta},$ we will show that%

\[
Kf(x)=\lim_{\varepsilon\rightarrow0}\int_{\varepsilon<d(x,y)}K(x,y)f(y)d\mu
(y)
\]
exists $\mu-a.e.$ \ Assume $\varepsilon<1,$ we can write

$\bigskip$%
\[
Kf(x)=\lim_{\varepsilon\rightarrow0}\int_{\varepsilon<d(x,y)<1}K(x,y) \left[
f(y)-f(x)\right]  d\mu(y)+
\]

\[
f(x)\lim_{\varepsilon\rightarrow0}\int_{\varepsilon<d(x,y)<1}K(x,y)d\mu
(y)+\int_{1\leq d(x,y)}K(x,y)f(y)d\mu(y).
\]
Note that the first integral converges absolutely, the limit of the second
term exists by hypothesis and finally last integral converges absolutely
because the integrand is bounded. \ Furthermore, we have $\left\Vert
Kf\right\Vert _{\infty}\leq c\left\Vert f\right\Vert _{Lip_{\beta}}.$

We will estimate now $Kf(x_{1})-Kf(x_{2})$ for $x_{1},x_{2}$ two points for
which $Kf(x)$ exists. \ This part of the proof is very similar to the same
part in Theorem 2. We write%

\[
Kf(x_{1})-Kf(x_{2})=\lim_{\varepsilon\rightarrow0}\int_{\varepsilon
<d(x_{1},y)}K(x_{1},y)f(y)d\mu(y)-\lim_{\varepsilon\rightarrow0}%
\int_{\varepsilon<d(x_{2},y)}K(x_{2},y)f(y)d\mu(y)
\]
\[
=\lim_{\varepsilon\rightarrow0}\int_{\varepsilon<d(x_{1},y)}K(x_{1},y) \left[
f(y)-f(x_{1})\right]  d\mu(y)+f(x_{1})\lim_{\varepsilon\rightarrow0}%
\int_{\varepsilon<d(x_{1},y)}K(x_{1},y)d\mu(y)-
\]

\[
\lim_{\varepsilon\rightarrow0}\int_{\varepsilon<d(x_{2},y)}K(x_{2},y) \left[
f(y)-f(x_{1})\right]  d\mu(y)-f(x_{1})\lim_{\varepsilon\rightarrow0}%
\int_{\varepsilon<d(x_{2},y)}K(x_{2},y)d\mu(y)=
\]

\[
\lim_{\varepsilon\rightarrow0}\int_{\varepsilon<d(x_{1},y)}K(x_{1},y) \left[
f(y)-f(x_{1})\right]  d\mu(y)+\lim_{\varepsilon\rightarrow0}\int
_{\varepsilon<d(x_{2},y)}K(x_{2},y)\left[  f(y)-f(x_{1})\right]  d\mu(y)+
\]

\[
f(x_{1})\left[  K1(x_{1})-K1(x_{2})\right]
\]

Observe now that the last term can be estimated using the hypothesis and we have%

\[
\left\vert f(x_{1})\left[  K1(x_{1})-K1(x_{2})\right]  \right\vert \leq
c\left\Vert f\right\Vert _{\infty}d^{\beta}(x_{1},x_{2}).
\]

To estimate the first two terms, let $r=d(x_{1},x_{2}),$ and $\varepsilon<r, $
we rewrite them as follows:%

\[
\lim_{\varepsilon\rightarrow0}\int_{\varepsilon<d(x_{1},y)}K(x_{1},y) \left[
f(y)-f(x_{1})\right]  d\mu(y)+\lim_{\varepsilon\rightarrow0}\int
_{\varepsilon<d(x_{2},y)}K(x_{2},y)\left[  f(y)-f(x_{1})\right]  d\mu(y)=
\]

\begin{align*}
&  \lim_{\varepsilon\rightarrow0}\int_{\varepsilon<d(x_{1},y)<3r}K(x_{1},y)
\left[  f(y)-f(x_{1})\right]  d\mu(y)+\\
&  \lim_{\varepsilon\rightarrow0}\int_{\left\{  y:\varepsilon<d(x_{2}%
,y)\right\}  \cap\left\{  y:d(x_{1},y)<3r\right\}  }K(x_{2},y)\left[
f(y)-f(x_{1})\right]  d\mu(y)+
\end{align*}

\[
\lim_{\varepsilon\rightarrow0}\int_{3r<d(x_{1},y)}\left[  f(y)-f(x_{1}%
)\right]  \left[  K(x_{1},y)-K(x_{2},y)\right]  d\mu(y)=H_{1}+H_{2}+H_{3}
\]
The absolute value of $H_{3}$ can be estimated as follows,%

\[
\left\vert H_{3}\right\vert \leq\left\vert f\right\vert _{\beta}d^{\gamma
}(x_{1},x_{2})\int_{3r<d(x_{1},y)}\frac{d^{\beta}(x_{1},y)}{d^{n+\gamma}%
(x_{1},y)}d\mu(y)\leq c\left\vert f\right\vert _{\beta}d^{\beta}(x_{1},x_{2})
\]

For $\left\vert H_{1}\right\vert $we have%

\[
\left\vert H_{1}\right\vert \leq\left\vert f\right\vert _{\beta}\int
_{d(x_{1},y)<3r}\frac{C_{1}}{d^{n-\beta}(x_{1},y)}d\mu(y)\leq c\left\vert
f\right\vert _{\beta}d^{\beta}(x_{1},x_{2})
\]

Finally to estimate $H_{2}$ we write%

\[
\lim_{\varepsilon\rightarrow0}\int_{\left\{  y:\varepsilon<d(x_{2},y)\right\}
\cap\left\{  y:d(x_{1},y)<3r\right\}  }K(x_{2},y)\left[  f(y)-f(x_{1})\right]
d\mu(y)=
\]

\[
\lim_{\varepsilon\rightarrow0}\int_{\left\{  y:\varepsilon<d(x_{2},y)\right\}
\cap\left\{  y:d(x_{1},y)<3r\right\}  }K(x_{2},y)\left[  f(y)-f(x_{2})\right]
d\mu(y)
\]

\[
+\left[  f(x_{2})-f(x_{1})\right]  \lim_{\varepsilon\rightarrow0}%
\int_{\left\{  y:\varepsilon<d(x_{2},y)\right\}  \cap\left\{  y:d(x_{1}%
,y)<3r\right\}  }K(x_{2},y)d\mu(y)=J_{1}+J_{2}
\]

For the first term we have%

\[
\left\vert J_{1}\right\vert \leq\int_{d(x_{2},y)<4r}\frac{c\left\Vert
f\right\Vert _{(\beta)}}{d^{n-\beta}(x_{2},y)}d\mu(y)\leq c\left\vert
f\right\vert _{\beta}d^{\beta}(x_{1},x_{2})
\]

To estimate the second $J_{2}$ consider first%

\[
\lim_{\varepsilon\rightarrow0}\int_{\left\{  y:\varepsilon<d(x_{2},y)\right\}
\cap\left\{  y:d(x_{1},y)<3r\right\}  }K(x_{2},y)d\mu(y)=\lim_{\varepsilon
\rightarrow0}\int_{\varepsilon<d(x_{2},y)<2r}K(x_{2},y)d\mu(y)+
\]

\[
\int_{\left\{  y:d(x_{1},y)<3r\right\}  \backslash\left\{  y:d(x_{2}%
,y)<2r\right\}  }K(x_{2},y)d\mu(y)
\]

Observe now that%

\[
\left\vert \lim_{\varepsilon\rightarrow0}\int_{\varepsilon<d(x_{2}%
,y)<2r}K(x_{2},y)d\mu(y)\right\vert \leq C_{3}
\]
and using part 3 of the lemma we get%

\[
\left\vert \int_{\left\{  y:d(x_{1},y)<3r\right\}  \backslash\left\{
y:d(x_{2},y)<2r\right\}  }K(x_{2},y)d\mu(y)\right\vert \leq\int_{\left\{
y:2r<d(x_{2},y)<4r\right\}  }\left\vert K(x_{2},y)\right\vert d\mu(y)\leq c
\]
therefore%

\[
\left\vert J_{2}\right\vert \leq c\left\vert f\right\vert _{\beta}d^{\beta
}(x_{1},x_{2})
\]
collecting the estimates we have:%

\[
\left\vert Kf(x_{1})-Kf(x_{2})\right\vert \leq c\left\Vert f\right\Vert
_{_{Lip_{\beta}}}d^{\beta}(x_{1},x_{2})
\]
and finally%

\[
\left\Vert Kf\right\Vert _{Lip_{\beta}}\leq c\left\Vert f\right\Vert
_{Lip_{\beta}}
\]
\bigskip This concludes the proof of Theorem 3.

\bigskip

\noindent\textbf{Proof of Theorem 4}

\qquad

\qquad We will prove the theorem for $D_{\alpha}(x,y)=$\bigskip$\frac
{1}{d^{n+\alpha}(x,y)},$the general case is identical. Note that the proof is
also valid for $\mu(X)=\infty.$

We will estimate first $\sup(D_{\alpha}f)$ for $f\in\Lambda_{\beta}.$We use
part 2 of the Lemma to write%

\[
\left\vert D^{\alpha}f(x)\right\vert \leq\int_{d(x,y)\leq1}\frac{\left\vert
f(y)-f(x)\right\vert }{d^{n+\alpha}(x,y)}d\mu(y)+c\sup(f)\leq\left\vert
f\right\vert _{\beta}\int_{d(x,y)\leq1}\frac{1}{d^{n+\alpha-\beta}(x,y)}%
d\mu(y)+c\sup(f)
\]
Since $0<\alpha<\beta\leq1,$ we use part 1 of the lemma to estimate the
integral and we obtain that $D^{\alpha}f(x)$ is well defined everywhere and%

\[
\sup(D^{\alpha}f)\leq c\left\Vert f\right\Vert _{\Lambda_{\beta}}.
\]

To estimate $\left\vert D^{\alpha}f\right\vert _{\alpha}$, we consider
$r=d(x,y)$ and write%

\begin{align*}
D^{\alpha}f(x_{1})-D^{\alpha}f(x_{2})  &  =\int_{d(x_{1},y)\leq2r}%
\frac{f(y)-f(x_{1})}{d^{n+\alpha}(x_{1},y)}d\mu(y)-\int_{d(x_{1},y)\leq
2r}\frac{f(y)-f(x_{2})}{d^{n+\alpha}(x_{2},y)}d\mu(y)+\\
&  +\int_{d(x_{1},y)>2r}\left[  f(y)-f(x_{1})\right]  \left[  \frac
{1}{d^{n+\alpha}(x_{1},y)}-\frac{1}{d^{n+\alpha}(x_{2},y)}\right]  d\mu(y)-\\
&  \int_{d(x_{1},y)>2r}\frac{f(x_{1})-f(x_{2})}{d^{n+\alpha}(x_{2},y)}d\mu(y)
\end{align*}

Using part1 of the Lemma and the fact that $f$ is in $\Lambda_{\beta}$ we can
obtain that each of the first two terms converges absolutely and is bounded by
$c\left\vert f\right\vert _{\beta}d^{\beta-\alpha}(x_{1},x_{2}).$ Using part 2
of the Lemma we can also obtain that the fourth term converges absolutely and
is bounded by $c\left\vert f\right\vert _{\beta}d^{\beta-\alpha}(x_{1},x_{2})$.

To estimate the third term observe first that for $2d(x_{1},x_{2})\leq
d(x_{1},y),$%

\[
\left\vert \frac{1}{d^{n+\alpha}(x_{1},y)}-\frac{1}{d^{n+\alpha}(x_{2}%
,y)}\right\vert \leq\sup_{\theta}\left\vert (-n-\alpha)(\theta d(x_{1}%
,y)+(1-\theta)(d(x_{2},y))^{-n-\alpha-1}\right\vert \text{ }.
\]

\[
\left\vert d(x_{1},y)-d(x_{2},y)\right\vert \leq c\frac{d(x_{1},x_{2}%
)}{d^{n+\alpha+1}(x_{1},y)}.
\]

Therefore using this estimate, the fact that $f\epsilon\Lambda_{\beta}$ and
the part 2 of Lemma we obtain that the third term converges absolutely and is
less than or equal to $c\left\vert f\right\vert _{\beta}d^{\beta-\alpha}%
(x_{1},x_{2})$ and consequently $\ \ \ \ \ \ \left\vert D^{\alpha}\right\vert
_{(\beta-\alpha)}\leq c\left\vert f\right\vert _{\beta}$.

Finally combining the two estimates we get $\left\Vert D^{\alpha}f\right\Vert
_{\Lambda_{\beta-\alpha}}\leq c\left\Vert f\right\Vert _{\Lambda_{\beta}}.$

\bigskip

To extend Theorem 1 and Theorem 3, the fractional integrals and singular
integrals have to be redefined so they converge for $d(x,y)>1.$The operator's
norm in each result will depend on the normalization.We will denote with ' the
normalizations. Let $x_{o}\in X$ be a fixed point for whuich (S4) is valid and define:%

\[
L_{\alpha}^{^{\prime}}f(x)=\int\left[  L_{\alpha}(x,y)-L_{\alpha}%
(x_{o},y)\right]  f(y)d\mu(y)
\]

\[
K^{\prime}f(x)=\lim_{\varepsilon\rightarrow0}\int_{\epsilon<d(x,y)}\left[
K(x,y)-K(x_{o},y)\right]  f(y)d\mu(y)
\]

\bigskip

\bigskip

\noindent\textbf{Applications}

\bigskip

In this section we will illustrate some applications of the theorems. I am
indebted to Joaquim Bruna for pointing out to me the Theorem of Mark Krein and
to Joan Verdera \bigskip for several generous discussions on applications 1
and 2.

1. The purpose of this application is to obtain boundedness in $L^{2}$ of some
singular integrals in the context of non-doubling measure metric spaces of
finite measure. Following $\left[  T\right]  ,$ a singular integral associated
to $\mu$ is said to be bounded in $L^{2}$ when there is a constant $C$ such
that $\left\Vert K_{\varepsilon}f\right\Vert _{L^{2}}\leq C\left\Vert
f\right\Vert _{L^{2}},$ for all $\varepsilon>0,$ where $K_{\varepsilon
}f(x)=\int_{d(x,y)>\varepsilon}K(x,y)f(y)d\mu(y)$. We will apply Theorem 2 and
the following Theorem of Mark Krein (see$\left[  \text{FMM}\right]  $ for its
proof and application to the classical case, and [W] for the case of spaces of
homogeneous type).

\noindent

M. Krein's Theorem:

Let $H$ be a real or complex Hilbert space with inner product $(.,.)$ and norm
$\left\Vert .\right\Vert _{H}.$ \ Let $D\subset H$ be a Banach space dense in
$H$ and such that $\left\Vert x\right\Vert _{H}\leq C\left\Vert x\right\Vert
_{D}$ for $x\in D.$ Let $A$ and $B$ be two linear operator such that
$\left\Vert Ax\right\Vert _{D}\leq C_{A}\left\Vert x\right\Vert _{D}$,
$\left\Vert Bx\right\Vert _{D}\leq C_{B}\left\Vert x\right\Vert _{D}$ $,x\in
D$\ \ and $(Ax,y)=(x,By)$ for all $x,y\in D$. \ Then $\left\Vert Ax\right\Vert
_{H}\leq(C_{A}C_{B})^{\frac{1}{2}}\left\Vert x\right\Vert _{H}$, $\left\Vert
Bx\right\Vert _{H}\leq(C_{A}C_{B})^{\frac{1}{2}}\left\Vert x\right\Vert _{H}$
$,x\in D,$ and both extend to bounded operator on $H.$

In our application, we will consider $H=L^{2}$ \ and $D=\Lambda_{\beta}$.
Since $X$ has finite measure we clearly have $\left\Vert f\right\Vert _{L^{2}%
}\leq\mu(X)^{\frac{1}{2}}\left\Vert f\right\Vert _{\Lambda_{\beta}}$, but we
need the extra assumption $\Lambda_{\beta}$ dense in $L^{2}.$ Let now $K(x,y)$
be a standard singular integral kernel and $K^{\ast}(x,y)=\overline{K(y,x)}$ .
Assume that $K^{\ast}(x,y)$ also satisfies (S2). Let $A=T_{\varepsilon}$ and
$B=T_{\varepsilon}^{\ast}$ the corresponding smooth truncations. If
$\left\Vert T_{\varepsilon}1\right\Vert _{\Lambda_{\beta}}\leq C^{\prime}$ and
$\left\Vert T_{\varepsilon}^{\ast}1\right\Vert _{\Lambda_{\beta}}\leq
C^{\prime\prime}$for all $\varepsilon>0,$ then by Theorem 2 and Krein's
Theorem there is $C$ such that $\left\Vert T_{\varepsilon}f\right\Vert
_{L^{2}}\leq C\left\Vert f\right\Vert _{L^{2}}$, $f\in\Lambda_{\beta}$, for
all $\varepsilon>0.$ Consequently$\left\Vert K_{\varepsilon}f\right\Vert
_{L^{2}}\leq C\left\Vert f\right\Vert _{L^{2}}$, $f\in\Lambda_{\beta}$ for all
$\varepsilon>0$, and it extends to a bounded operator in $L^{2}$,same
conclusion for $K^{\ast}$. In addition, Nazarov, Treil, and Volberg have
extended the classical result of Calderon-Zygmund on the boundedness in
$L^{p},$ $1<p<\infty,$of singular integrals bounded in $L^{2},$ to
non-doubling separable measure metric spaces, see $\left[  \text{NTV}\right]
.$

2. The second application has appeared in $\left[  \text{MOV}\right]  .$ In
this paper the authors need to study the boundedness properties of the
Restricted Beurling Transform$,$ $B_{\Omega}f=B(f\chi_{\Omega}),$on
$Lip_{\varepsilon}(\Omega)$ where $\Omega$ is a bounded domain in $R^{n}$ with
boundary of class $C^{1+\varepsilon}$, $0<\varepsilon<1,.$ Mateu, Orobitg and
Verdera prove the following more general result: \ "Let $\Omega$ be a bounded
domain with boundary of class $C^{1+\varepsilon}$,$0<\varepsilon<1$, and let
$T$ be an even smooth homogeneous Calderon-Zygmund operator. \ Then
$T_{\Omega}$ maps $Lip_{\varepsilon}\left(  \Omega\right)  $ into
$Lip_{\varepsilon}(\Omega)$ and also $Lip_{\varepsilon}\left(  \Omega\right)
$ into $Lip_{\varepsilon}(\Omega^{c})$". \ Their proof, which is non-trivial,
consists in showing that condition (S3) and part a) of Theorem 3 above are met.

3. The third application is related to M. Riesz Fractional Calculus associated
to non-doubling measures. Applying Theorem 1 and Theorem 4 we can obtain that
the composition of a Riesz fractional integral $I_{\alpha}f(x)=\int\frac
{1}{d^{n-\alpha}(x,y)}f(y)d\mu(y)$ and a fractional \ derivative $D^{\alpha
}f(x)=\int\frac{\left[  f(y)-f(x)\right]  }{d^{n+\alpha}(x,y)}d\mu(y)$ of the
same order $D^{\alpha}I_{\alpha}$, as well as its transpose $I_{\alpha
}D^{\alpha}$, are bounded on $\Lambda_{\beta},$ when $I_{\alpha}1$ $\in
\Lambda_{\alpha+\beta},\alpha+\beta<1.$\ In addition, it was shown in $\left[
\text{G}\right]  $ that these composition are singular integral operators
associated to $\mu$\pagebreak

\bigskip

\bigskip

e-mail of the author:

aegatto@depaul.edu

\end{document}